\documentclass[3p]{elsarticle}

\usepackage{hyperref, amssymb, amsmath, amsthm, subfig,color}
\usepackage[font=scriptsize]{caption}
\newtheorem{theorem}{Theorem}[section]
\newtheorem{lemma}[theorem]{Lemma}

\newtheorem{corollary}[theorem]{Corollary}

\theoremstyle{definition}

\theoremstyle{remark}

\usepackage{pgfplots}
\pgfplotsset{compat=1.18}

\journal{}

\makeatletter
\def\ps@pprintTitle{%
 \let\@oddhead\@empty
 \let\@evenhead\@empty
 \def\@oddfoot{}%
 \let\@evenfoot\@oddfoot}
\makeatother

\begin{document}

\begin{frontmatter}

\title{New upper and lower bounds on the smallest singular values of\\nonsingular lower triangular $(0,1)$-matrices}

\author{Vesa Kaarnioja}
\ead{vesa.kaarnioja@lut.fi}
\address{School of Engineering Sciences, LUT University, P.O.~Box 20, 53851 Lappeenranta, Finland}
\author{Andr\'e-Alexander Zepernick\corref{cor1}}\cortext[cor1]{Corresponding author}
\ead{a.zepernick@fu-berlin.de}
\address{Freie Universit\"at Berlin, Fachbereich Mathematik und Informatik, Arnimallee 6, 14195 Berlin, Germany}

\begin{abstract}
Let $K_n$ denote the set of all nonsingular $n\times n$ lower triangular $(0,1)$-matrices. Hong and Loewy (2004) introduced the number sequence
$$
c_n=\min\{\lambda\mid\lambda~\text{is an eigenvalue of}~XX^{\rm T},~X\in K_n\},\quad n\in\mathbb Z_+.
$$
There have been a number of attempts in the literature to obtain bounds on the numbers $c_n$ by Mattila (2015), Alt{\i}n{\i}\c{s}{\i}k et al.~(2016), Kaarnioja (2021), Loewy (2021), and Alt{\i}n{\i}\c{s}{\i}k (2021). In this paper, improved upper and lower bounds are derived for the numbers $c_n$. By considering the characteristic polynomial corresponding to the matrix $Z_n$ satisfying $c_n=\|Z_n\|_2^{-1}$, it is shown that the second largest eigenvalue of $Z_n$ is bounded from above by $\frac45$ leading to an improved upper bound on $c_n$. On the other hand, Samuelson's inequality applied to the roots of the characteristic polynomial of $Z_n$ yields an improved lower bound. Numerical experiments demonstrate the quality of the new bounds.
\end{abstract}


\end{frontmatter}

\section{Introduction}
Let $K_n$ denote the set of all nonsingular lower triangular $(0,1)$-matrices. The set $K_2$ contains the matrices
$$
\begin{pmatrix}
1&0\\0&1
\end{pmatrix},~\begin{pmatrix}
1&0\\1&1
\end{pmatrix},
$$
while the set $K_3$ consists of the matrices\begin{align*}
&\begin{pmatrix}
1&0&0\\
0&1&0\\
0&0&1
\end{pmatrix},~\begin{pmatrix}
1&0&0\\
1&1&0\\
0&0&1
\end{pmatrix},~\begin{pmatrix}
1&0&0\\
0&1&0\\
1&0&1
\end{pmatrix},~\begin{pmatrix}
1&0&0\\
0&1&0\\
0&1&1
\end{pmatrix},\\
&\begin{pmatrix}
1&0&0\\
1&1&0\\
1&0&1
\end{pmatrix},~\begin{pmatrix}
1&0&0\\
1&1&0\\
0&1&1
\end{pmatrix},~\begin{pmatrix}
1&0&0\\
0&1&0\\
1&1&1
\end{pmatrix},~\begin{pmatrix}
1&0&0\\
1&1&0\\
1&1&1
\end{pmatrix}.
\end{align*}
In fact, the cardinality of $K_n$ is $2^{n(n-1)/2}$ for all $n\in\mathbb{Z}_+$. 

Hong and Loewy~\cite{hongloewy} introduced the numbers
$$
c_n=\min\{\lambda\mid\lambda~\text{is an eigenvalue of}~XX^{\rm T},~X\in K_n\},\quad n\in\mathbb Z_+,
$$
in order to obtain lower bounds for the eigenvalues of greatest common divisor (GCD) matrices defined on a set of arbitrary integers. However, Hong and Loewy did not give bounds on the numbers $c_n$ in~\cite{hongloewy}. In recent years, a number of papers have been published about developing explicit bounds on the numbers $c_n$, see~\cite{mattila,altinisik,kaarnioja2021,loewy2021,altinisik21}. Besides GCD matrices, the numbers $c_n$ also appear in eigenvalue bounds for least common multiple (LCM) matrices and the more general classes of meet and join matrices~\cite{ihm}.

The numbers $c_n$ are related to the singular values of matrices belonging to $K_n$. Letting $\sigma_{\min}(X)$ denote the smallest singular value of matrix $X$, there holds
$$
\min_{X\in K_n}\sigma_{\min}(X)=\sqrt{c_n}.
$$
Moreover, these numbers can be used to derive bounds on the smallest eigenvalues of lattice-theoretic matrices. As an example, let $(P,\preceq,\wedge,\hat 0)$ be a locally finite meet semilattice, where $\preceq$ is a partial ordering on the set $P$, $\wedge$ is the \emph{meet} (\emph{greatest lower bound}) of two elements in $P$, and $\hat 0\in P$ is the least element in $P$ such that $\hat 0\preceq x$ for all $x\in P$. If $S=\{x_1,\ldots,x_n\}\subset P$ is a lower closed set such that $x_i\preceq x_j$ only if $i\leq j$ and $f\!:P\to \mathbb R$ is a function, then we define the $n\times n$ meet matrix $A$ by setting $A_{i,j}=f(x_i\wedge x_j)$ for $i,j\in\{1,\ldots,n\}$ and the function
$$
J_{P,f}(x)=\sum_{\hat 0\preceq z\preceq x}f(z)\mu_P(z,x)\quad\text{for all}~x\in P,
$$
where $\mu_P$ denotes the M\"obius function of $P$ defined inductively by setting
$$
\mu_P(x,y)=\begin{cases}
1&\text{if}~x=y,\\
-\sum_{x\preceq z\prec y}\mu_P(x,z)&\text{if}~x\prec y,\\
0&\text{otherwise}.
\end{cases}
$$
If $J_{P,f}(x)>0$ for all $x\in P$, then
\begin{align}\label{eq:prelim}
\lambda_{\min}(A)\geq c_n\min_{x\in S}J_{P,f}(x),
\end{align}
where $\lambda_{\min}(A)$ denotes the smallest eigenvalue of matrix $A$. The inequality~\eqref{eq:prelim} was derived for meet matrices associated with incidence functions in~\cite{ihm}, but a lower bound of the form~\eqref{eq:prelim} already appeared in~\cite{hongloewy} within the context of the divisor lattice $(\mathbb Z_+,|,\gcd)$ and function $f(x)=x$.

The spectrum of meet and join matrices has been considered by many authors. Upper and lower bounds on the extremal eigenvalues of Smith's matrix $(\gcd(i,j))_{i,j=1}^n$ were derived by Balatoni~\cite{balatoni}. Positive definiteness of GCD matrices defined over an arbitrary set of distinct positive integers was proved by Beslin and Ligh~\cite{beslinligh}. Bourque and Ligh determined conditions which ensure that GCD matrices composed with arithmetic functions are positive definite~\cite{bourqueligh}. Hong studied the behavior of the largest eigenvalues of GCD matrices associated with certain arithmetic functions~\cite{hong2008}, while the behavior of the smallest eigenvalue was considered by Hong and Loewy~\cite{hongloewy11}. Hong and Lee investigated the asymptotic behavior of the eigenvalues of reciprocal power LCM matrices~\cite{hongenochlee2008}. Merikoski improved the lower bound on the smallest eigenvalue of Smith's matrix~\cite{merikoski}, eigenvalue bounds for ``mixed'' power GCD and LCM matrices were derived by Mattila and Haukkanen~\cite{haukkanenmattila12}, improved eigenvalue bounds for GCD and LCM matrices were derived by Alt{\i}n{\i}\c{s}{\i}k and B\"{u}y\"{u}kk\"{o}se~\cite{altinisik16}, and Ilmonen derived bounds on the eigenvalues of meet hypermatrices~\cite{ilmonen}. Haukkanen et al.~gave a lattice-theoretic interpretation for the inertia of LCM matrices~\cite{inertia}.

Computing the numbers $c_n$ directly from definition becomes intractable even for moderate values of $n$. However,  Alt{\i}n{\i}\c{s}{\i}k et al.~\cite{altinisik} showed that
\begin{align}
c_n=\|Z_n\|_2^{-1},\label{eq:spectral}
\end{align}
where the $n\times n$ matrix $Z_n$ is defined elementwise by setting
\begin{align}\label{eq:Zdef}
(Z_n)_{i,j}=\begin{cases}
1+\sum_{k=i+1}^nF_{k-i}^2&\text{if }i=j,\\
(-1)^{j-i}(F_{j-i}+\sum_{k=j+1}^nF_{k-i}F_{k-j})&\text{if }i<j,\\
(-1)^{i-j}(F_{i-j}+\sum_{k=i+1}^nF_{k-i}F_{k-j})&\text{if }i>j,
\end{cases}
\end{align}
for $i,j\in\{1,\ldots,n\}$ with the sequence of Fibonacci numbers defined by the recurrence $F_1=F_2=1$ and $F_k=F_{k-1}+F_{k-2}$ for $k\geq 3$. This was the basis for the lower bound derived on the numbers $c_n$ in~\cite{kaarnioja2021}. Let $\|\,\cdot\,\|_2$ denote the matrix spectral norm and let $\|\,\cdot\,\|_{\rm F}$ denote the Frobenius norm.
Since $\|\,\cdot\,\|_2\leq \|\,\cdot\,\|_{\rm F}$, the lower bound obtained in~\cite{kaarnioja2021} was based on the inequality
\begin{align}
c_n\geq \|Z_n\|_{\rm F}^{-1},\label{eq:lower}
\end{align}
where the main contribution was to verify the explicit identity
\begin{align}
\|Z_n\|_{\rm F}\!=\!\sqrt{\frac{1}{25}\varphi^{-4n}\!+\!\frac{3+(-1)^n}{25}\varphi^{-2n}\!-\!\frac{2}{5\sqrt 5}n\varphi^{-2n}\!+\!\frac{13(-1)^n\!-\!33}{50}\!+\!n\!+\!\frac{3\!+\!(-1)^n}{25}\varphi^{2n}\!+\!\frac{2}{5\sqrt 5}n\varphi^{2n}\!+\!\frac{1}{25}\varphi^{4n}},\label{eq:Z_nFrob}
\end{align}
with $\varphi=\frac{1+\sqrt 5}{2}$ denoting the golden ratio. The quality of the lower bound~\eqref{eq:lower} was demonstrated by numerical experiments. Shortly thereafter, Loewy~\cite{loewy2021} derived the bounds
\begin{align}
\frac{2}{2F_n^2+2n-1+(-1)^n}\leq c_n\leq \frac{2}{2F_n^2+1+(-1)^n},\quad n\in\mathbb Z_+,\label{eq:loewybound}
\end{align}
proving that $c_n\sim 5\varphi^{-2n}$ as $n\to\infty$.

The goal of this paper is to give a new upper bound on the numbers $c_n$, which improves Loewy's bounds~\eqref{eq:loewybound}. The lower bound previously derived in~\cite{kaarnioja2021} is improved in the present work by conducting a more refined analysis on the spectrum of $Z_n$.

This paper is structured as follows. An explicit expression for the characteristic polynomial of $Z_n$ is derived in Section~\ref{sec:2}, leading us to conclude that its second largest eigenvalue can be bounded from above by $\frac45$. This property is used to obtain a new upper bound on the numbers $c_n$ in Section~\ref{sec:3}, while a novel lower bound on the numbers $c_n$ is derived in Section~\ref{sec:4} by bounding the roots of the characteristic polynomial of $Z_n$. Numerical experiments validating the quality of these new bounds are presented in Section~\ref{sec:5} and some conclusions summarizing our findings are presented in Section~\ref{sec:6}.

\section{An upper bound on the second largest eigenvalue of $Z_n$}\label{sec:2}

Let $Z_n$ be defined by~\eqref{eq:Zdef}. It will turn out to be necessary to control the behavior of the characteristic polynomial $p_n(\lambda)=\det(\lambda I-Z_n)$ at $\lambda=\frac45$. We will first derive a general formula for the characteristic polynomial.
\begin{lemma}\label{lemma:lemma21}
    For all $n\in\mathbb{Z}_+$, there holds
    \begin{align}
    p_n(\lambda)&=\frac12 (\lambda-1)^{n-1}+(\lambda-1)^n-\frac{\lambda (1-\lambda)^n}{8-18\lambda+10\lambda^2}\notag
    \\&\quad+\frac{2^{-n}(-2+3\lambda-\sqrt\lambda \sqrt{5\lambda-4})^n+2^{-n}(-2+3\lambda +\sqrt\lambda \sqrt{5\lambda-4})^n}{4-5\lambda}.\label{eq:induc}
    \end{align}
\end{lemma}
\proof  Denote $H_n=\lambda I-Z_n$ and consider for the moment the case $n\geq 5$. The determinant is
\begin{align*}&\det H_n=\\&
\resizebox{\textwidth}{!}{$
\begin{vmatrix}
(\lambda-1)-\sum_{k=2}^n F_{k-1}F_{k-1} & F_1 + \sum_{k=3}^n F_{k-2}F_{k-1} & -F_2-\sum_{k=4}^n F_{k-3}F_{k-1} & F_3+\sum_{k=5}^n F_{k-4}F_{k-1} & -F_4-\sum_{k=6}^n F_{k-5}F_{k-1} & \cdots & (-1)^n F_{n-1} \\
F_1+\sum_{k=3}^n F_{k-1}F_{k-2} & (\lambda-1)-\sum_{k=3}^n F_{k-2}F_{k-2} & F_1+\sum_{k=4}^n F_{k-3}F_{k-2} & -F_2-\sum_{k=5}^n F_{k-4}F_{k-2} & F_3+\sum_{k=6}^n F_{k-5}F_{k-2} & \cdots & (-1)^{n-1}F_{n-2} \\
-F_2-\sum_{k=4}^n F_{k-1}F_{k-3} & F_1+\sum_{k=4}^n F_{k-2}F_{k-3} & (\lambda-1)-\sum_{k=4}^n F_{k-3}F_{k-3} & F_1+\sum_{k=5}^n F_{k-4}F_{k-3} & -F_2-\sum_{k=6}^n F_{k-5}F_{k-3} & \cdots & (-1)^{n-2}F_{n-3} \\
F_3+\sum_{k=5}^n F_{k-1}F_{k-4} & -F_2-\sum_{k=5}^n F_{k-2}F_{k-4} & F_1+\sum_{k=5}^n F_{k-3}F_{k-4} & (\lambda-1)-\sum_{k=5}^n F_{k-4}F_{k-4} & F_1+\sum_{k=6}^n F_{k-5}F_{k-4} & \cdots & (-1)^{n-3}F_{n-4} \\
-F_4-\sum_{k=6}^n F_{k-1}F_{k-5} & F_3+\sum_{k=6}^n F_{k-2}F_{k-5} & -F_2-\sum_{k=6}^n F_{k-3}F_{k-5} & F_1+\sum_{k=6}^n F_{k-4}F_{k-5} & (\lambda-1)-\sum_{k=6}^n F_{k-5}F_{k-5} & \cdots & (-1)^{n-4}F_{n-5} \\
\vdots & \vdots & \vdots & \vdots & \vdots & \ddots & \vdots
\end{vmatrix}.
$}
\end{align*}
By performing the determinant-preserving elementary row and column operations
\begin{align*}
&{\rm row}_1\rightarrow {\rm row}_1+{\rm row}_2-{\rm row}_3,\\
&{\rm col}_1\rightarrow {\rm col}_1+{\rm col}_2-{\rm col}_3,
\end{align*}
this determinant can be written as
\begin{align*}&\det H_n=\\&
\resizebox{\textwidth}{!}{$
\begin{vmatrix}
3\lambda-2 &\lambda-1 &1-\lambda &0&0& \cdots & 0\\
\lambda-1& (\lambda-1)-\sum_{k=3}^n F_{k-2}F_{k-2} & F_1+\sum_{k=4}^n F_{k-3}F_{k-2} & -F_2-\sum_{k=5}^n F_{k-4}F_{k-2} & F_3+\sum_{k=6}^n F_{k-5}F_{k-2} & \cdots & (-1)^{n-1}F_{n-2} \\
1-\lambda & F_1+\sum_{k=4}^n F_{k-2}F_{k-3} & (\lambda-1)-\sum_{k=4}^n F_{k-3}F_{k-3} & F_1+\sum_{k=5}^n F_{k-4}F_{k-3} & -F_2-\sum_{k=6}^n F_{k-5}F_{k-3} & \cdots & (-1)^{n-2}F_{n-3} \\
0 & -F_2-\sum_{k=5}^n F_{k-2}F_{k-4} & F_1+\sum_{k=5}^n F_{k-3}F_{k-4} & (\lambda-1)-\sum_{k=5}^n F_{k-4}F_{k-4} & F_1+\sum_{k=6}^n F_{k-5}F_{k-4} & \cdots & (-1)^{n-3}F_{n-4} \\
0& F_3+\sum_{k=6}^n F_{k-2}F_{k-5} & -F_2-\sum_{k=6}^n F_{k-3}F_{k-5} & F_1+\sum_{k=6}^n F_{k-4}F_{k-5} & (\lambda-1)-\sum_{k=6}^n F_{k-5}F_{k-5} & \cdots & (-1)^{n-4}F_{n-5} \\
\vdots & \vdots & \vdots & \vdots & \vdots & \ddots & \vdots
\end{vmatrix}.
$}
\end{align*}
In particular, the new matrix contains $H_{n-1}$ as a principal submatrix. Further, performing the elementary row and column operations
\begin{align*}
&{\rm row}_2\rightarrow {\rm row}_2+{\rm row}_3-{\rm row}_4,\\
&{\rm col}_2\rightarrow {\rm col}_2+{\rm col}_3-{\rm col}_4,\\
&{\rm row}_3\rightarrow {\rm row}_3+{\rm row}_4-{\rm row}_5,\\
&{\rm col}_3\rightarrow {\rm col}_3+{\rm col}_4-{\rm col}_5,
\end{align*}
results in
\begin{align*}
&\det H_n=\\
&\resizebox{\textwidth}{!}{$\begin{vmatrix}
3\lambda-2 &0 &1-\lambda &0&0& \cdots & 0\\
0& 3\lambda-2& 0 &1-\lambda & 0& \cdots & 0\\
1-\lambda & 0 & 3\lambda-2& \lambda-1 & 1-\lambda& \cdots & 0 \\
0 & 1-\lambda &\lambda-1& (\lambda-1)-\sum_{k=5}^n F_{k-4}F_{k-4} & F_1+\sum_{k=6}^n F_{k-5}F_{k-4} & \cdots & (-1)^{n-3}F_{n-4} \\
0& 0& 1-\lambda & F_1+\sum_{k=6}^n F_{k-4}F_{k-5} & (\lambda-1)-\sum_{k=6}^n F_{k-5}F_{k-5} & \cdots & (-1)^{n-4}F_{n-5} \\
\vdots & \vdots & \vdots & \vdots & \vdots & \ddots & \vdots
\end{vmatrix},$}
\end{align*}
which contains $H_{n-3}$ as a principal submatrix.
By carrying out the Laplace cofactor expansion across the first column, we obtain
\begin{align*}
\det H_n=&\,(3\lambda-2)\det H_{n-1}\\
&+(1-\lambda)\resizebox{.78\textwidth}{!}{$ \begin{vmatrix}
0 & 1-\lambda & 0 & 0 & \cdots & 0\\
3\lambda-2 & 0 & 1-\lambda & 0 & \cdots & 0\\
1-\lambda & \lambda-1&(\lambda-1)-\sum_{k=5}^n F_{k-4}F_{k-4} & F_1+\sum_{k=6}^n F_{k-5}F_{k-4} &\cdots & (-1)^{n-3}F_{n-4} \\
0 & 1-\lambda&F_1+\sum_{k=6}^n F_{k-4}F_{k-5} & (\lambda-1)-\sum_{k=6}^n F_{k-5}F_{k-5} & \cdots & (-1)^{n-4}F_{n-5} \\
\vdots & \vdots & \vdots& \vdots & \ddots & \vdots
\end{vmatrix}$.}
\end{align*}
Next we carry out the Laplace cofactor expansion across the first row, which yields
\begin{align}
p_n(\lambda)&=\det H_n\notag\\
&=(3\lambda-2)\det H_{n-1}\notag\\
&\quad -(1-\lambda)^2\resizebox{.78\textwidth}{!}{$\begin{vmatrix}
3\lambda-2 & 1-\lambda & 0 & \cdots & 0\\
1-\lambda & (\lambda-1)-\sum_{k=5}^n F_{k-4}F_{k-4} & F_1+\sum_{k=6}^n F_{k-5}F_{k-4}&\cdots & (-1)^{n-3}F_{n-4} \\
0 & F_1+\sum_{k=6}^n F_{k-4}F_{k-5} & (\lambda-1)-\sum_{k=6}^n F_{k-5}F_{k-5}&\cdots & (-1)^{n-4}F_{n-5} \\
\vdots & \vdots &\vdots& \ddots & \vdots
\end{vmatrix}$}\notag\\
&=(3\lambda-2)\det H_{n-1}-(1-\lambda)^2\big((3\lambda-2)\det H_{n-3}+(1-\lambda)^2\det H_{n-4}\big)\notag\\
&=(3\lambda-2)\det H_{n-1}-(1-\lambda)^2(3\lambda-2)\det H_{n-3}+(1-\lambda)^4\det H_{n-4}\notag\\
&=(3\lambda-2)(\det H_{n-1}-(1-\lambda)^2\det H_{n-3})+(1-\lambda)^4\det H_{n-4}\notag\\
&=(3\lambda-2)(p_{n-1}(\lambda)-(1-\lambda)^2p_{n-3}(\lambda))+(1-\lambda)^4p_{n-4}(\lambda),\label{detdecu}
\end{align}
where the third equality is a consequence of performing the Laplace cofactor expansion across the first row of the remaining determinant.
With the recurrence formula~\eqref{detdecu} on hand, it remains to use induction to prove the identity~\eqref{eq:induc}. For $n=1,\ldots,4$, using basic operations, it is straightforward to compute the characteristic polynomials
$$
p_1(\lambda)=\lambda-1,\quad p_2(\lambda)=\lambda^2-3\lambda+1,\quad p_3(\lambda)=\lambda^3-6\lambda^2+5\lambda-1,\quad p_4(\lambda)=\lambda^4-13\lambda^3+18\lambda^2-8\lambda+1,
$$
which agree with the formula~\eqref{eq:induc} by simplifying the expressions, which one obtains for the different values of $n$. Fix the integer $n\geq 5$ and assume that the claim has already been proved for all positive integers $\leq n-1$. By the derived recurrence formula we have
\begin{align*}
  p_n(\lambda)&= (3\lambda-2)(p_{n-1}(\lambda)-(1-\lambda)^2p_{n-3}(\lambda))+(1-\lambda)^4p_{n-4}(\lambda).
\end{align*} 
Exploiting the induction hypothesis yields
\begin{align*}
    p_{n-1}(\lambda) &= \frac{1}{2} (\lambda - 1)^{n-2} + (\lambda - 1)^{n-1} - \frac{\lambda (1 - \lambda)^{n-1}}{8 - 18\lambda + 10\lambda^2}\\
   &\quad + \frac{2^{-n+1} (-2 + 3\lambda - \sqrt{\lambda} \sqrt{5\lambda - 4})^{n-1} + 2^{-n+1} (-2 + 3\lambda + \sqrt{\lambda} \sqrt{5\lambda - 4})^{n-1}}{4 - 5\lambda},\\
   p_{n-3}(\lambda) &= \frac{1}{2} (\lambda - 1)^{n-4} + (\lambda - 1)^{n-3} - \frac{\lambda (1 - \lambda)^{n-3}}{8 - 18\lambda + 10\lambda^2}\\
   &\quad+ \frac{2^{-n+3} (-2 + 3\lambda - \sqrt{\lambda} \sqrt{5\lambda - 4})^{n-3} + 2^{-n+3} (-2 + 3\lambda + \sqrt{\lambda} \sqrt{5\lambda - 4})^{n-3}}{4 - 5\lambda},\quad \text{and}\\
   p_{n-4}(\lambda) &= \frac{1}{2} (\lambda - 1)^{n-5} + (\lambda - 1)^{n-4} - \frac{\lambda (1 - \lambda)^{n-4}}{8 - 18\lambda + 10\lambda^2}\\
   &\quad+ \frac{2^{-n+4} (-2 + 3\lambda - \sqrt{\lambda} \sqrt{5\lambda - 4})^{n-4} + 2^{-n+4} (-2 + 3\lambda + \sqrt{\lambda} \sqrt{5\lambda - 4})^{n-4}}{4 - 5\lambda}.
\end{align*}
We need to show that plugging in these expressions into the recurrence relation yields the desired formula. Therefore we look at the terms involving powers of $(\lambda-1)$. Noting that $(1-\lambda)^2=(\lambda-1)^2$ as well as $(1-\lambda)^4=(\lambda-1)^4$ we find
\begin{align*}
  &(3\lambda-2)\left(\frac{1}{2}(\lambda-1)^{n-2}+(\lambda-1)^{n-1}\right)-(3\lambda-2)(1-\lambda)^2\left(\frac{1}{2}(\lambda-1)^{n-4}+(\lambda-1)^{n-3}\right)\\&+(1-\lambda)^4\left(\frac{1}{2}(\lambda-1)^{n-5}+(\lambda-1)^{n-4}\right)\\
  &=\frac{(3\lambda-2)(\lambda-1)^{n-2}}{2}+(3\lambda-2)(\lambda-1)^{n-1}-\frac{(3\lambda-2)(\lambda-1)^{n-2}}{2}-(3\lambda-2)(\lambda-1)^{n-1}\\
  &\quad +\frac{(\lambda-1)^{n-1}}{2}+(\lambda-1)^n\\
  &=\frac{1}{2}(\lambda-1)^{n-1}+(\lambda-1)^n.
\end{align*}
Furthermore, we observe
\begin{align*}
    &-(3\lambda-2)\frac{\lambda(1-\lambda)^{n-1}}{8-18\lambda+10\lambda^2}+(3\lambda-2)(1-\lambda)^2\frac{\lambda(1-\lambda)^{n-3}}{8-18\lambda+10\lambda^2}-(1-\lambda)^4\frac{\lambda(1-\lambda)^{n-4}}{8-18\lambda+10\lambda^2}\\
    &=-(3\lambda-2)\frac{\lambda(1-\lambda)^{n-1}}{8-18\lambda+10\lambda^2}+(3\lambda-2)\frac{\lambda(1-\lambda)^{n-1}}{8-18\lambda+10\lambda^2}-\frac{\lambda(1-\lambda)^{n}}{8-18\lambda+10\lambda^2}\\
    &=-\frac{\lambda(1-\lambda)^{n}}{8-18\lambda+10\lambda^2}.
\end{align*}
Finally, we have
\begin{align*}
    & (3\lambda - 2) \cdot \frac{2^{-n+1} \left( -2 + 3\lambda - \sqrt{\lambda} \sqrt{5\lambda - 4} \right)^{n-1} + 2^{-n+1} \left( -2 + 3\lambda + \sqrt{\lambda} \sqrt{5\lambda - 4} \right)^{n-1}}{4 - 5\lambda} \\
    & \quad - (3\lambda - 2)(1 - \lambda)^2 \cdot \frac{2^{-n+3} \left( -2 + 3\lambda - \sqrt{\lambda} \sqrt{5\lambda - 4} \right)^{n-3} + 2^{-n+3} \left( -2 + 3\lambda + \sqrt{\lambda} \sqrt{5\lambda - 4} \right)^{n-3}}{4 - 5\lambda} \\
    & \quad + (1 - \lambda)^4 \cdot \frac{2^{-n+4} \left( -2 + 3\lambda - \sqrt{\lambda} \sqrt{5\lambda - 4} \right)^{n-4} + 2^{-n+4} \left( -2 + 3\lambda + \sqrt{\lambda} \sqrt{5\lambda - 4} \right)^{n-4}}{4 - 5\lambda}\\
    &=\frac{2^{-n+1}(3\lambda-2)(A^{n-1}+B^{n-1})-2^{-n+3}(3\lambda-2)(1-\lambda)^2(A^{n-3}+B^{n-3})+2^{-n+4}(1-\lambda)^4(A^{n-4}+B^{n-4})}{4-5\lambda},
\end{align*}
where we define
\begin{align*}
    A&:=-2+3\lambda -\sqrt{\lambda}\sqrt{5\lambda-4},\\
    B&:=-2+3\lambda +\sqrt{\lambda}\sqrt{5\lambda-4}.
\end{align*}
It remains to show that
\begin{align*}
 2(3\lambda-2)(A^{n-1}+B^{n-1})-8(3\lambda-2)(1-\lambda)^2(A^{n-3}+B^{n-3})+16(1-\lambda)^4(A^{n-4}+B^{n-4} )=A^n+B^n.  
\end{align*}
Note that $A$ and $B$ fulfill the equations
\begin{align*}
    A+B=2 (-2+3\lambda)=6\lambda-4\enspace \text{and}\enspace A\cdot B=4 (1-\lambda)^2.
\end{align*}
Hence, by Vieta's theorem \cite[p.~102]{vieta} we conclude that $A$ and $B$ are the roots of the polynomial
\begin{align*}
    x^2-(6\lambda-4)x+4(1-\lambda)^2.
\end{align*}
Since this can be interpreted as the characteristic polynomial of a recurrence relation, it follows that $A$ and $B$ fulfill the recurrence formula
\begin{align*}A^n+B^n=(6\lambda-4)(A^{n-1}+B^{n-1})-4(1-\lambda)^2(A^{n-2}+B^{n-2}).
\end{align*}
Expanding the terms further yields
\begin{align*}
  A^n+B^n&= (6\lambda-4)(A^{n-1}+B^{n-1})-4(1-\lambda)^2(A^{n-2}+B^{n-2})\\
  &=(6\lambda-4)(A^{n-1}+B^{n-1})-4(1-\lambda)^2\big((6\lambda-4)(A^{n-3}+B^{n-3})-4(1-\lambda)^2(A^{n-4}+B^{n-4})\big)\\
  &=(6\lambda-4)(A^{n-1}+B^{n-1})-4(1-\lambda)^2(6\lambda-4)(A^{n-3}+B^{n-3})+16(1-\lambda)^4(A^{n-4}+B^{n-4})\\
  &=2(3\lambda-2)(A^{n-1}+B^{n-1})-8(3\lambda-2)(1-\lambda)^2(A^{n-3}+B^{n-3})+16(1-\lambda)^4(A^{n-4}+B^{n-4} ),
\end{align*}
which verifies the desired identity. Combining all terms from above proves the statement.
\endproof
The behavior of the characteristic polynomial $p_n(\lambda)$ at $\lambda=\frac45$ plays a special role in what follows.
\begin{corollary} For all $n\in\mathbb Z_+$, there holds
\begin{align}
p_n(\tfrac45)=\det (\tfrac45 I-Z_n)=-\frac12\cdot \frac{1}{5^n}(3(-1)^n+10n^2-5).\label{eq:induc2}
\end{align}
\end{corollary}

\proof Although plugging the value $\lambda=\frac45$ into the expression~\eqref{eq:induc} would yield the assertion, simplifying the resulting expression turns out to be rather technical. We present in the following a simpler argument.

From~\eqref{detdecu}, we obtain in the special case $\lambda=\frac45$ the recurrence
\begin{align}
p_n(\tfrac45)&=\frac25p_{n-1}(\tfrac45)-\frac{2}{125}p_{n-3}(\tfrac45)+\frac{1}{625}p_{n-4}(\tfrac45),\quad n\geq 5.\label{eq:recu2}
\end{align}
Similarly to the proof of Lemma~\ref{lemma:lemma21}, we use induction to prove the identity~\eqref{eq:induc2}. For $n=1,\ldots,4$, we obtain 
$$
p_1(\tfrac45)=-\frac15,\quad p_2(\tfrac45)=-\frac{19}{25},\quad p_3(\tfrac45)=-\frac{41}{125},\quad p_4(\tfrac45)=-\frac{79}{625},
$$
which coincide with~\eqref{eq:induc2}. We fix the integer $n\geq 5$ and assume that the claim has already been proved for all integers $\leq n-1$. Then the recurrence~\eqref{eq:recu2} yields
\begin{align*}
p_n(\tfrac45)&=-\frac{1}{5^n}(3(-1)^{n-1}+10(n-1)^2-5) + \frac{1}{5^{n}}(3(-1)^{n-3}+10(n-3)^2-5)\\
&\quad -\frac{1}{2}\cdot \frac{1}{5^{n}}(3(-1)^{n-4}+10(n-4)^2-5))\\
&=-\frac12\cdot\frac{1}{5^n}(3(-1)^n+10(2(n-1)^2-2(n-3)^2+(n-4)^2)-5\big)
\end{align*}
and the claim follows by expanding
$$
2(n-1)^2-2(n-3)^2+(n-4)^2=n^2,
$$
completing the induction.\endproof

\begin{lemma}\label{lemma:eigbnd}
Let $n\geq 2$ be an integer. The eigenvalues $\lambda_1^{(n)}\geq \lambda_2^{(n)}\geq \cdots\geq \lambda_n^{(n)}>0$ of $Z_n$ satisfy $\lambda_2^{(n)}<\frac45$ and $\lambda_1^{(n)}>1$.
\end{lemma}
\proof The proof is carried out by induction with respect to~$n$. When $n=2$, the eigenvalues of $Z_2\in\mathbb R^{2\times 2}$ are given by
$$
\lambda_1^{(2)}=\frac{3+\sqrt 5}{2}>1\quad\text{and}\quad \lambda_2^{(2)}=\frac{3-\sqrt 5}{2}<\frac45,
$$
which resolves the basis of the induction.

Next, suppose that $\lambda_2^{(n-1)}<\frac45$ and $\lambda_1^{(n-1)}>1$ for some $n$. Since
$$
Z_n=\begin{pmatrix}*&*\\ *&Z_{n-1}\end{pmatrix},
$$
i.e., the matrix $Z_n$ contains $Z_{n-1}$ as a principal submatrix with eigenvalues $\lambda_1^{(n-1)}\geq\lambda_2^{(n-1)}\geq\cdots\geq \lambda_{n-1}^{(n-1)}>0 $, it is a direct consequence of the Cauchy interlacing theorem (cf., e.g., \cite[Theorem~4.3.17]{cauchy}) that
\begin{align}\label{eq:cauchyinter}
\lambda_k^{(n)}\leq \lambda_{k}^{(n-1)}\leq \lambda_{k+1}^{(n)}\quad\text{for all}~k\in\{1,\ldots,n-1\}.
\end{align}
In particular, it follows that $\lambda_1^{(n)}\geq \lambda_1^{(n-1)}>1$ and $0<\lambda_n^{(n)}\leq \cdots\leq \lambda_{2}^{(n)}\leq \lambda_2^{(n-1)}<\frac45$. Consider the characteristic polynomial $$p_n(\lambda)=\det(\lambda I-Z_n)=(\lambda-\lambda_1^{(n)})(\lambda - \lambda_2^{(n)})(\lambda - \lambda_3^{(n)})\cdots (\lambda-\lambda_n^{(n)}).$$
Plugging in the value $\lambda=\frac45$ reveals that
$$p_n(\tfrac45)=\underset{<0}{\underbrace{(\tfrac45-\lambda_1^{(n)})}}(\tfrac45 - \lambda_2^{(n)})\underset{>0}{\underbrace{(\tfrac45 - \lambda_3^{(n)})\cdots (\tfrac45-\lambda_n^{(n)})}}.$$
However, we already know that
$$
p_n(\tfrac45)=\det (\tfrac45 I-Z_n)<0\quad\text{for all}~n\in\mathbb Z_+,
$$
which implies that $\frac45-\lambda_2^{(n)}>0$ and hence $\lambda_2^{(n)}<\frac45$.\endproof

\emph{Remark.} The inequality $\lambda_1^{(n)}>1$ for the largest eigenvalue with $n\geq 2$ was already proved in~\cite{loewy2021}, but this is reproved in Lemma~\ref{lemma:eigbnd} for completeness. Furthermore, we note that it is a consequence of~\eqref{eq:cauchyinter} that $\lambda_2^{(2)}\leq \lambda_2^{(3)}\leq \lambda_2^{(4)}\leq\cdots\leq \frac45$.

\section{An improved upper bound on Hong and Loewy's numbers $c_n$}\label{sec:3}

Let $n\geq 2$ be an integer and let $\lambda_1^{(n)}\geq \lambda_2^{(n)}\geq \cdots\geq \lambda_n^{(n)}>0$ denote the eigenvalues of the positive definite matrix $Z_n$. By Lemma~\ref{lemma:eigbnd}, there holds $0<\lambda_n^{(n)}\leq \cdots\leq \lambda_2^{(n)}<\frac45$. Hence $\|Z_n\|_{\rm F}^2=(\lambda_1^{(n)})^2+\cdots+(\lambda_n^{(n)})^2<\|Z_n\|_2^2+\frac{16}{25}(n-1)$. Similarly to~\cite[Theorem~2.1]{loewy2021} and~\cite[Theorem~2]{altinisik21} we obtain by a sequence of simple algebraic manipulations that
\begin{align*}
\|Z_n\|_{\rm F}^2<\|Z_n\|_2^2+\frac{16}{25}(n-1)\quad &\Leftrightarrow \quad \frac{1}{\|Z_n\|_2^2}<\frac{1}{\|Z_n\|_{\rm F}^2}+\frac{\frac{16}{25}(n-1)}{\|Z_n\|_{\rm F}^2\|Z_n\|_2^2}\\
&\Leftrightarrow\quad \frac{1}{\|Z_n\|_2^2}\bigg(1+\frac{\frac{16}{25}(1-n)}{\|Z_n\|_{\rm F}^2}\bigg)<\frac{1}{\|Z_n\|_{\rm F}^2}\\
&\Leftrightarrow\quad \frac{1}{\|Z_n\|_2}<\frac{1}{\sqrt{\|Z_n\|_{\rm F}^2+\frac{16}{25}(1-n)}}.
\end{align*}
Plugging in the identities~\eqref{eq:spectral} and~\eqref{eq:Z_nFrob} yields the following result. 
\begin{theorem}\label{thm:main}
There holds for all integers $n\geq 2$ that
\begin{align}
c_n<\frac{1}{\sqrt{\frac{1}{25}\varphi^{-4n}+\frac{3+(-1)^n}{25}\varphi^{-2n}-\frac{2}{5\sqrt 5}n\varphi^{-2n}+\frac{13(-1)^n-33}{50}+1+\frac{3+(-1)^n}{25}\varphi^{2n}+\frac{2}{5\sqrt 5}n\varphi^{2n}+\frac{1}{25}\varphi^{4n}}}.\label{eq:ub}
\end{align}
\end{theorem}
\section{A new lower bound on Hong and Loewy's numbers $c_n$}\label{sec:4}
The lower bound on $c_n$ derived in~\cite{kaarnioja2021} can be improved by obtaining bounds on the roots of the characteristic polynomial $p_n(\lambda)= \det (\lambda I-Z_n)$ instead of exploiting the matrix norm inequality~\eqref{eq:lower}.
\begin{theorem}\label{thm:newlower}
Let $n\geq 2$ be an integer. Then
\begin{align}
&c_n\notag\\
&\!\geq \!\frac{1}{1\!\hspace*{-.02cm}+\hspace*{-.02cm}\!\frac{(-1)^n\!-\!5}{10n}\!\hspace*{-.02cm}+\hspace*{-.02cm}\!\frac{\varphi^{2n}\!+\!\varphi^{-2n}}{5n}\!\hspace*{-.02cm}+\hspace*{-.02cm}\!\frac{n\!-\!1}{5n}\!\sqrt{\!\frac{4(-1)^n}{n\!-\!1}\!\hspace*{-.02cm}+\hspace*{-.02cm}\!\frac{2\sqrt 5 n^2(\varphi^{2n}\!-\!\varphi^{-2n})}{n\!-\!1}\!\hspace*{-.02cm}+\hspace*{-.02cm}\! \frac{3(-1)^n\!+\!17}{2}\!\hspace*{-.02cm}+\hspace*{-.02cm}\!((-1)^n\!\hspace*{-.02cm}-\hspace*{-.02cm}\!7\!\hspace*{-.02cm}-\hspace*{-.02cm}\!\frac{2}{n\!-\!1})(\varphi^{2n}\!\hspace*{-.02cm}+\hspace*{-.02cm}\!\varphi^{-2n})\!\hspace*{-.02cm}+\hspace*{-.02cm}\!\varphi^{4n}\!\hspace*{-.02cm}+\hspace*{-.02cm}\!\varphi^{-4n}}}.\label{eq:lb}
\end{align}
\end{theorem}
\proof Let $p_n(x)=a_nx^n+a_{n-1}x^{n-1}+\cdots a_1x+a_0$, $a_n\neq 0$, be a polynomial with $n$ real roots. Samuelson's inequality~\cite{samuelson} implies that the roots of $p_n$ are located in the interval
$$
\bigg[-\frac{a_{n-1}}{na_n}-\frac{n-1}{na_n}\sqrt{a_{n-1}^2-\frac{2n}{n-1}a_na_{n-2}},-\frac{a_{n-1}}{na_n}+\frac{n-1}{na_n}\sqrt{a_{n-1}^2-\frac{2n}{n-1}a_na_{n-2}}\bigg].
$$
In particular, it is sufficient to know the coefficients of the three leading terms of the polynomial $p_n$ to obtain an upper bound on its largest root.

If $p_n(\lambda)=\det (\lambda I-Z_n)$, then it is an immediate consequence of Newton's identities~\cite{newton} that
$$
a_n=1,\quad a_{n-1}=-{\rm tr}(Z_n),\quad \text{and}\quad a_{n-2}=\frac{1}{2}({\rm tr}(Z_n)^2-\|Z_n\|_{\rm F}^2).
$$
By plugging in the value~\eqref{eq:Z_nFrob} for $\|Z_n\|_{\rm F}$ and the value
$$
{\rm tr}(Z_n)=n+F_n^2-\frac{1-(-1)^n}{2}=n-\frac12+\frac15\varphi^{2n}+\frac15\varphi^{-2n}+\frac{1}{10}(-1)^n,
$$
one obtains
\begin{align*}
&c_n\geq \frac{1}{-\frac{a_{n-1}}{na_n}+\frac{n-1}{na_n}\sqrt{a_{n-1}^2-\frac{2n}{n-1}a_na_{n-2}}}\\
&\!=\!\frac{1}{1\!\hspace*{-.02cm}+\hspace*{-.02cm}\!\frac{(-1)^n\!-\!5}{10n}\!\hspace*{-.02cm}+\hspace*{-.02cm}\!\frac{\varphi^{2n}\!+\!\varphi^{-2n}}{5n}\!\hspace*{-.02cm}+\hspace*{-.02cm}\!\frac{n\!-\!1}{5n}\!\sqrt{\!\frac{4(-1)^n}{n\!-\!1}\!\hspace*{-.02cm}+\hspace*{-.02cm}\!\frac{2\sqrt 5 n^2(\varphi^{2n}\!-\!\varphi^{-2n})}{n\!-\!1}\!\hspace*{-.02cm}+\hspace*{-.02cm}\! \frac{3(-1)^n\!+\!17}{2}\!\hspace*{-.02cm}+\hspace*{-.02cm}\!((-1)^n\!\hspace*{-.02cm}-\hspace*{-.02cm}\!7\!\hspace*{-.02cm}-\hspace*{-.02cm}\!\frac{2}{n\!-\!1})(\varphi^{2n}\!\hspace*{-.02cm}+\hspace*{-.02cm}\!\varphi^{-2n})\!\hspace*{-.02cm}+\hspace*{-.02cm}\!\varphi^{4n}\!\hspace*{-.02cm}+\hspace*{-.02cm}\!\varphi^{-4n}}}
\end{align*}
as desired.\endproof
\section{Numerical experiments}\label{sec:5}

We assess the quality of the lower bound presented in Theorem~\ref{thm:newlower} and the upper bound derived in Theorem~\ref{thm:main}. In Table~\ref{table:table}, the values of $c_n$ have been tabulated against the lower bounds given in Theorem~\ref{thm:newlower} vis-\`a-vis the lower bound~\eqref{eq:lower}--\eqref{eq:Z_nFrob} obtained in~\cite{kaarnioja2021} for $n\in\{2,\ldots,10\}$. We also display the value of the upper bound in~Theorem~\ref{thm:main} in Table~\ref{table:table}. The numerical values suggest that the lower bound of Theorem~\ref{thm:newlower} is an improvement over~\eqref{eq:lower}--\eqref{eq:Z_nFrob}. Indeed, the lower and upper bounds become better for $n>10$. Note that the lower bound given by Theorem~\ref{thm:newlower} coincides with the true value of $c_n$ for $n=2$. 
\begin{table}[!t]
\centering
\begin{tabular}{c|c|cc|c}
$n$&$c_n$&$\overset{\text{\normalsize Lower bound}}{\text{in Theorem~\ref{thm:newlower}}}$&$\overset{\text{\normalsize Lower bound~\eqref{eq:lower}--\eqref{eq:Z_nFrob}}}{\text{derived in~\cite{kaarnioja2021}}}$&$\overset{\text{\normalsize Upper bound}}{\text{in Theorem~\ref{thm:main}}}$\\
\hline$2$&$0.3819660113$&$0.3819660113$&$0.3779644730$&$0.4082482905$\\
$3$&$0.1980622642$&$0.1978219619$&$0.1961161351$&$0.2041241452$\\
$4$&$0.0870031120$&$0.0869565217$&$0.0867109970$&$0.0877058019$\\
$5$&$0.0370683347$&$0.0370629486$&$0.0370370370$&$0.0371390676$\\
$6$&$0.0148275852$&$0.0148271154$&$0.0148249863$&$0.0148331386$\\
$7$&$0.0058169987$&$0.0058169617$&$0.0058168052$&$0.0058173957$\\
$8$&$0.0022453455$&$0.0022453429$&$0.0022453322$&$0.0022453719$\\
$9$&$0.0008622031$&$0.0008622030$&$0.0008622023$&$0.0008622048$\\
$10$&$0.0003300037$&$0.0003300037$&$0.0003300036$&$0.0003300038$
\end{tabular}
\caption{Tabulated values of the constant $c_n$, the lower bound in Theorem~\ref{thm:newlower}, the previous lower bound~\eqref{eq:lower}--\eqref{eq:Z_nFrob} derived in~\cite{kaarnioja2021}, and the upper bound in Theorem~\ref{thm:main} for $n\in\{2,\ldots,10\}$.}\label{table:table}
\end{table}
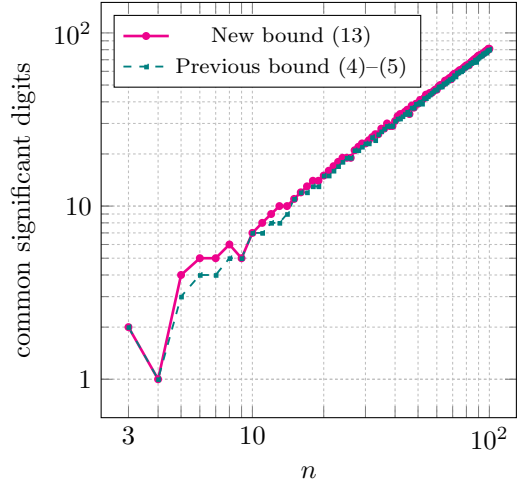
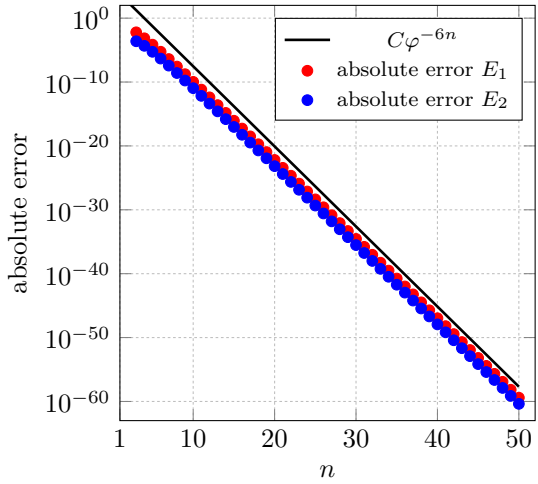
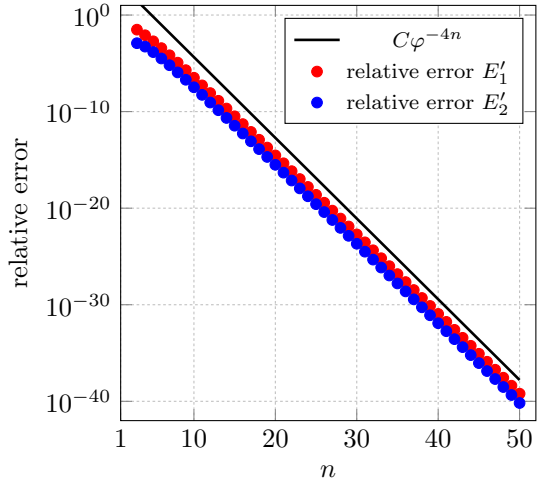
\begin{figure}[!t]\centering
\subfloat[]{\label{fig:Xa}\begin{tikzpicture}
    \begin{loglogaxis}[
        width=.43\textwidth,
        height=.43\textwidth,
        xlabel={$n$},
        ylabel={common significant digits},
        legend pos=north west,
        grid=both,
        grid style={dash pattern=on 1pt off 1pt on 1pt off 1pt},
        xmin=2.3, xmax=130,
        ymin=0.6, ymax=150,
        xtick={2,3,4,5,6,7,8,9,10,20,30,40,50,60,70,80,90,100},
        xticklabels={,3,,,,,,,10,,,,,,,,,$10^2$},
        ytick={1,10,100},
        yticklabels={1, 10, $10^2$}
    ]
    \addplot[
        color=magenta,
        mark=* ,
        mark size=1pt,
        line width=1pt
    ]
    coordinates {
        (3,2) (4,1) (5,4) (6,5) (7,5) (8,6) (9,5) (10,7) (11,8) (12,9)
        (13,10) (14,10) (15,11) (16,12) (17,13) (18,14) (19,14) (20,15)
        (21,16) (22,17) (23,18) (24,19) (25,19) (26,19) (27,21) (28,22)
        (29,23) (30,23) (31,24) (32,25) (33,26) (34,26) (35,28) (36,28)
        (37,30) (38,29) (39,29) (40,31) (41,33) (42,34) (43,34) (44,35)
        (45,36) (46,34) (47,38) (48,37) (49,39) (50,39) (51,41) (52,41)
        (53,42) (54,44) (55,44) (56,45) (57,45) (58,46) (59,47) (60,47)
        (61,49) (62,50) (63,50) (64,51) (65,53) (66,53) (67,54) (68,55)
        (69,54) (70,57) (71,58) (72,58) (73,59) (74,60) (75,61) (76,61)
        (77,62) (78,63) (79,64) (80,65) (81,65) (82,65) (83,67) (84,68)
        (85,69) (86,70) (87,71) (88,71) (89,73) (90,74) (91,74) (92,75)
        (93,75) (94,76) (95,77) (96,78) (97,79) (98,79) (99,81) (100,81)
    };
    \addlegendentry{\text{\footnotesize New bound~\eqref{eq:lb}}}
    \addplot[
        color=teal,
        mark=square* ,
        mark size=0.7pt,
        line width=0.7pt,
        dashed
    ]
    coordinates {
        (3,2) (4,1) (5,3) (6,4) (7,4) (8,5) (9,5) (10,7) (11,7) (12,8)
        (13,8) (14,9) (15,11) (16,12) (17,12) (18,13) (19,13) (20,15)
        (21,15) (22,16) (23,17) (24,18) (25,19) (26,19) (27,21) (28,21)
        (29,22) (30,23) (31,23) (32,25) (33,24) (34,26) (35,27) (36,28)
        (37,29) (38,29) (39,29) (40,31) (41,32) (42,32) (43,33) (44,34)
        (45,35) (46,34) (47,37) (48,37) (49,39) (50,39) (51,40) (52,39)
        (53,42) (54,42) (55,43) (56,44) (57,45) (58,46) (59,47) (60,47)
        (61,49) (62,49) (63,49) (64,51) (65,51) (66,52) (67,53) (68,54)
        (69,54) (70,55) (71,57) (72,56) (73,59) (74,59) (75,60) (76,60)
        (77,61) (78,63) (79,63) (80,64) (81,65) (82,65) (83,67) (84,67)
        (85,68) (86,68) (87,70) (88,68) (89,71) (90,72) (91,73) (92,74)
        (93,74) (94,76) (95,76) (96,77) (97,77) (98,79) (99,80) (100,80)
    };
    \addlegendentry{\text{\footnotesize Previous bound~\eqref{eq:lower}--\eqref{eq:Z_nFrob}}}
    \end{loglogaxis}
\end{tikzpicture}}
\\
\subfloat[]{\label{fig:Xb}\begin{tikzpicture}
    \begin{semilogyaxis}[
        width=.43\textwidth,
        height=.43\textwidth,
        xlabel={$n$},
        ylabel={absolute error},
        legend pos=north east,
        grid=both,
        grid style={dash pattern=on 1pt off 1pt on 1pt off 1pt},
        xmin=1, xmax=52,
        ymin=1e-63, ymax=1e2,
        xtick={1, 10, 20, 30, 40,50},
        ytick={1e-60, 1e-50, 1e-40, 1e-30, 1e-20, 1e-10, 1e0}
    ]
    \addplot[
        color=black,
        domain=2:50,
        samples=100,
        line width=1pt
    ] {100000*((1 + sqrt(5))/2)^(-6*x)};
    \addplot[
        color=red,
        mark=*,
        only marks
    ] coordinates {
        (3, 0.00606188)
        (4, 0.00070269)
        (5, 0.0000707329)
        (6, 5.5534e-6)
        (7, 3.96971e-7)
        (8, 2.63677e-8)
        (9, 1.69236e-9)
        (10, 1.06108e-10)
        (11, 6.55707e-12)
        (12, 4.0085e-13)
        (13, 2.43037e-14)
        (14, 1.46363e-15)
        (15, 8.76465e-17)
        (16, 5.2231e-18)
        (17, 3.09947e-19)
        (18, 1.83244e-20)
        (19, 1.07979e-21)
        (20, 6.34406e-23)
        (21, 3.71743e-24)
        (22, 2.17308e-25)
        (23, 1.26754e-26)
        (24, 7.37875e-28)
        (25, 4.28758e-29)
        (26, 2.48721e-30)
        (27, 1.44059e-31)
        (28, 8.33194e-33)
        (29, 4.81254e-34)
        (30, 2.77629e-35)
        (31, 1.59975e-36)
        (32, 9.20813e-38)
        (33, 5.29481e-39)
        (34, 3.0417e-40)
        (35, 1.74579e-41)
        (36, 1.00116e-42)
        (37, 5.73676e-44)
        (38, 3.28476e-45)
        (39, 1.87944e-46)
        (40, 1.07464e-47)
        (41, 6.14064e-49)
        (42, 3.50671e-50)
        (43, 2.0014e-51)
        (44, 1.14163e-52)
        (45, 6.5086e-54)
        (46, 3.70877e-55)
        (47, 2.11232e-56)
        (48, 1.20251e-57)
        (49, 6.84268e-59)
        (50, 3.89204e-60)
    };
    \addplot[
        color=blue,
        mark=*,
        only marks
    ] coordinates {
        (3, 0.000240302)
        (4, 0.0000465902)
        (5, 5.3861e-6)
        (6, 4.69806e-7)
        (7, 3.70588e-8)
        (8, 2.57663e-9)
        (9, 1.71444e-10)
        (10, 1.09335e-11)
        (11, 6.83988e-13)
        (12, 4.20997e-14)
        (13, 2.56538e-15)
        (14, 1.55027e-16)
        (15, 9.30886e-18)
        (16, 5.55946e-19)
        (17, 3.3051e-20)
        (18, 1.95706e-21)
        (19, 1.15479e-22)
        (20, 6.79282e-24)
        (21, 3.98461e-25)
        (22, 2.33147e-26)
        (23, 1.36109e-27)
        (24, 7.92948e-29)
        (25, 4.61083e-30)
        (26, 2.67645e-31)
        (27, 1.55111e-32)
        (28, 8.97602e-34)
        (29, 5.18715e-35)
        (30, 2.99378e-36)
        (31, 1.72582e-37)
        (32, 9.93775e-39)
        (33, 5.71649e-40)
        (34, 3.28508e-41)
        (35, 1.8861e-42)
        (36, 1.08196e-43)
        (37, 6.20153e-45)
        (38, 3.55184e-46)
        (39, 2.03279e-47)
        (40, 1.1626e-48)
        (41, 6.64481e-50)
        (42, 3.79546e-51)
        (43, 2.16665e-52)
        (44, 1.23614e-53)
        (45, 7.04872e-55)
        (46, 4.01726e-56)
        (47, 2.28842e-57)
        (48, 1.30298e-58)
        (49, 7.4155e-60)
        (50, 4.21849e-61)
    };
    \legend{\text{\footnotesize$C\varphi^{-6n}$},\text{\footnotesize absolute error $E_1$},\text{\footnotesize absolute error $E_2$}}
    \end{semilogyaxis}
\end{tikzpicture}}
\subfloat[]{\label{fig:Xc}\begin{tikzpicture}
    \begin{semilogyaxis}[
        width=.43\textwidth,
        height=.43\textwidth,
        xlabel={$n$},
        ylabel={relative error},
        legend pos=north east,
        grid=both,
        grid style={dash pattern=on 1pt off 1pt on 1pt off 1pt},
        xmin=1, xmax=52,
        ymin=1e-42, ymax=1e1,
        xtick={1, 10, 20, 30, 40,50},
        ytick={1e-60, 1e-50, 1e-40, 1e-30, 1e-20, 1e-10, 1e0}
    ]
    \addplot[
        color=black,
        domain=2:50,
        samples=100,
        line width=1pt
    ] {10000*((1 + sqrt(5))/2)^(-4*x)};
    \addplot[
        color=red,
        mark=*,
        only marks
    ] coordinates {
        (3, 0.0306059)
        (4, 0.00807661)
        (5, 0.00190818)
        (6, 0.000374532)
        (7, 0.0000682433)
        (8, 0.0000117433)
        (9, 1.96283e-6)
        (10, 3.21535e-7)
        (11, 5.19695e-8)
        (12, 8.3145e-9)
        (13, 1.31956e-9)
        (14, 2.08035e-10)
        (15, 3.26138e-11)
        (16, 5.08822e-12)
        (17, 7.90493e-13)
        (18, 1.22353e-13)
        (19, 1.88756e-14)
        (20, 2.90337e-15)
        (21, 4.45403e-16)
        (22, 6.8165e-17)
        (23, 1.04093e-17)
        (24, 1.58642e-18)
        (25, 2.41337e-19)
        (26, 3.66522e-20)
        (27, 5.5578e-21)
        (28, 8.41558e-22)
        (29, 1.27259e-22)
        (30, 1.922e-23)
        (31, 2.89946e-24)
        (32, 4.36929e-25)
        (33, 6.57756e-26)
        (34, 9.89249e-27)
        (35, 1.48647e-27)
        (36, 2.23174e-28)
        (37, 3.34797e-29)
        (38, 5.01873e-30)
        (39, 7.51788e-31)
        (40, 1.12539e-31)
        (41, 1.68357e-32)
        (42, 2.51705e-33)
        (43, 3.76098e-34)
        (44, 5.61653e-35)
        (45, 8.38311e-36)
        (46, 1.25061e-36)
        (47, 1.86478e-37)
        (48, 2.77929e-38)
        (49, 4.14043e-39)
        (50, 6.16555e-40)
    };
    \addplot[
        color=blue,
        mark=*,
        only marks
    ] coordinates {
        (3, 0.00121327)
        (4, 0.000535501)
        (5, 0.000145302)
        (6, 0.0000316846)
        (7, 6.37077e-6)
        (8, 1.14754e-6)
        (9, 1.98845e-7)
        (10, 3.31315e-8)
        (11, 5.42109e-9)
        (12, 8.73238e-10)
        (13, 1.39286e-10)
        (14, 2.20349e-11)
        (15, 3.46389e-12)
        (16, 5.4159e-13)
        (17, 8.42939e-14)
        (18, 1.30674e-14)
        (19, 2.01866e-15)
        (20, 3.10875e-16)
        (21, 4.77416e-17)
        (22, 7.31336e-18)
        (23, 1.11776e-18)
        (24, 1.70483e-19)
        (25, 2.59532e-20)
        (26, 3.94408e-21)
        (27, 5.98419e-22)
        (28, 9.06612e-23)
        (29, 1.37165e-23)
        (30, 2.07257e-24)
        (31, 3.12794e-25)
        (32, 4.71549e-26)
        (33, 7.10139e-27)
        (34, 1.06841e-27)
        (35, 1.60594e-28)
        (36, 2.41184e-29)
        (37, 3.61921e-30)
        (38, 5.42681e-31)
        (39, 8.13126e-32)
        (40, 1.21751e-32)
        (41, 1.82179e-33)
        (42, 2.72431e-34)
        (43, 4.07151e-35)
        (44, 6.08147e-36)
        (45, 9.07879e-37)
        (46, 1.35464e-37)
        (47, 2.02024e-38)
        (48, 3.01148e-39)
        (49, 4.48703e-40)
        (50, 6.68268e-41)
    };
    \legend{\text{\footnotesize$C\varphi^{-4n}$},\text{\footnotesize relative error $E_1'$},\text{\footnotesize relative error $E_2'$}}
    \end{semilogyaxis}
\end{tikzpicture}}
\caption{Comparisons of the lower bound~\eqref{eq:lb} and upper bound~\eqref{eq:ub} against the true value of $c_n$ and the previous lower bound~\eqref{eq:lower}--\eqref{eq:Z_nFrob} derived in~\cite{kaarnioja2021} for increasing $n$.}
\end{figure}

The number of common significant digits between $c_n$ and the lower bound in Theorem~\ref{thm:newlower} has been plotted in Figure~\ref{fig:Xa}. This is compared against the number of significant digits corresponding to the previous bound~\eqref{eq:lower}--\eqref{eq:Z_nFrob} derived in~\cite{kaarnioja2021}. We observe that the new lower bound is always at least as good as the previous lower bound, and both bounds become sharper as $n$ increases. In addition, we compute the differences between the lower bound~\eqref{eq:lb} against $c_n$ and the upper bound~\eqref{eq:ub} against $c_n$, respectively. We compute the absolute and relative errors
\begin{align*}
&E_1=(\text{upper bound}~\eqref{eq:ub})-c_n,\quad E_1'=\frac{E_1}{c_n},\\
&E_2=c_n-(\text{lower bound}~\eqref{eq:lb}),\quad\text{and}\quad E_2'=\frac{E_2}{c_n}. 
\end{align*}
These are illustrated in Figures~\ref{fig:Xb} and~\ref{fig:Xc}. Both the upper and lower difference tend to zero at an exponential rate, meaning that the bounds capture the correct asymptotic behavior of the numbers $c_n$. All numerical experiments were carried out using 150 digit precision computations in Wolfram Mathematica 12.3.

\section{Conclusions}\label{sec:6}
Theorem~\ref{thm:main} improves the upper bound on $c_n$ given in~\cite{loewy2021} considerably. Together with the new lower bound given in Theorem~\ref{thm:newlower}, these constitute improved bounds for the numbers $c_n$ over those previously considered in the literature. It remains an open challenge to the community whether even better bounds on $c_n$ can be obtained or if it is possible to produce an analytic expression for the numbers $c_n$ for all $n\in\mathbb Z_+$. However, since $c_n$ corresponds to a root of a polynomial of degree $n$, there is no good reason to expect that a closed form expression exists. In addition, numerical evidence suggests that the second largest eigenvalues $\lambda_2^{(n)}$ tend toward $\frac45$ very rapidly. Quantifying this convergence rate is deferred to future work.

\section*{Declaration of competing interest}
The authors declare that they have no known competing financial interests or personal relationships that could have appeared to influence the work reported in this paper.
\section*{Data availability}
No data was used for the research described in the article.

\bibliography{boolean_improvement}
\end{document}